\documentclass[10pt]{article}
\usepackage{lmodern}
\usepackage{amsmath}
\usepackage[T1]{fontenc}
\usepackage[utf8]{inputenc}
\usepackage{authblk}
\usepackage{amsfonts}
\usepackage{graphicx}
\usepackage{rotating}
\usepackage{amssymb}
\usepackage[english]{babel}
\usepackage{xcolor}
\usepackage{amsthm}
\usepackage{graphicx}
\usepackage{mathrsfs}
\usepackage{makecell}
\usepackage{microtype}
\usepackage{mathscinet}
\usepackage{array}
\usepackage{multirow}
\usepackage{enumerate}
\usepackage[cal=boondoxo,bb=ams]{mathalfa}
\usepackage{hyperref}
\hypersetup{hidelinks}
\usepackage{booktabs}

\newtheorem{defn}{Definition}[section]
\newtheorem{theorem}{Theorem}[section]
\newtheorem{prop}{Proposition}[section]
\newtheorem{lemma}{Lemma}[section]

\newtheorem{remark}{Remark}[section]

\newcommand{\ml}{\mathcal}
\newcommand{\mb}{\mathbb}
\DeclareMathOperator{\G}{\mathbb{G}}
\DeclareMathOperator{\lin}{lin}

\def\XXint#1#2#3{{\setbox0=\hbox{$#1{#2#3}{\int}$ }
		\vcenter{\hbox{$#2#3$ }}\kern-.6\wd0}}

\title{Blow-up and sharp lifespan estimates for a weakly coupled system of semilinear wave equations on a compact Lie group}
\author[1]{Wenhui Chen\thanks{Wenhui Chen (wenhui.chen.math@gmail.com)}}
\affil[1]{School of Mathematics and Information Science, Guangzhou University,\authorcr 510006 Guangzhou, P.R. China}

\author[2]{Alessandro Palmieri\thanks{Alessandro Palmieri (alessandro.palmieri@uniba.it)}}
\affil[2]{Department of Mathematics, University of Bari, 70125 Bari, Italy}

\date{}

\begin{document}
		\maketitle

		\begin{abstract}
			\medskip

In this paper, we investigate the blow-up in finite time and the corresponding lifespan estimates for a weakly coupled system of wave equations on a compact Lie group. In particular, we show how the Cauchy data and the presence of lower order terms affect the lifespan of local in-time solutions. 
			\\
			
			
\noindent\textbf{AMS Classification (2020)}  35B44, 35L05, 35L56, 35L76, 58J45
		\end{abstract}
\fontsize{12}{15}
\selectfont

\section{Introduction}\setcounter{equation}{0}\label{Section-Introduction}
\hspace{5mm}Let $\G$ be a compact Lie group (equipped with the uniquely determined normalized bi-invariant Riemannian metric) and let $\ml{L}$ be the Laplace-Beltrami operator on $\G$. 
In the present paper, we consider the following weakly coupled system of a semilinear wave equation and a semilinear 
damped Klein-Gordon equation:
\begin{align}\label{Eq-Main-Systems}
\begin{cases}
\partial_t^2u-\ml{L}u+b\,\partial_tu+m^2u=|v|^p,&x\in\G,\ t>0,\\
\partial_t^2v-\ml{L}v=|u|^q,&x\in\G,\ t>0,\\
(u,\partial_tu)(0,x)=\varepsilon(u_0,u_1)(x),&x\in\G,\\
(v,\partial_tv)(0,x)=\varepsilon(v_0,v_1)(x),&x\in\G,
\end{cases}
\end{align}
where $p,q>1$, $\varepsilon$ is a positive constant describing the size of Cauchy data, and the constant coefficients for the damping and mass terms satisfy $b>0$, $m^2\geqslant0$ as well as the further assumption 
\begin{align}\label{dominant damping assumption}
b^2\geqslant 4m^2.
\end{align}
In a nutshell, we can describe the condition \eqref{dominant damping assumption} as a constraint for the coefficients of the damping and mass terms in the $u$-equation. It guarantees the dominance of the damping term over the mass term: by this claim we mean that some properties of the solution to the linearized $u$-equation when \eqref{dominant damping assumption} holds are more similar to those for the damped wave equation (i.e. for $m^2=0$) than to those for the Klein-Gordon equation (i.e. for $b=0$).

The main goal of the present manuscript is to derive blow-up results and (sharp) lifespan estimates for the solutions to the weakly coupled system \eqref{Eq-Main-Systems} both for the case $b^2\geqslant 4m^2>0$ and for the case $b^2>4m^2=0$.

 In the following, whenever we refer to a nonlinearity of power-type, we always mean a non-negative power nonlinearity ($|u|^p$ or $|\partial_t u|^p$, for instance). We briefly review some literature related to the problem in \eqref{Eq-Main-Systems} in the Euclidean setting. Nowadays, it is well established that the critical curves for the weakly coupled system of semilinear wave equations (cf. \cite{ AgKuTa2000, De97, DeGeMi1997, DeMi1998, GeTaZh2006, IkSoWa2019, Ku05, KuTa03, KuTaWa2012}) and for the weakly coupled system of semilinear damped wave equations (cf. \cite{Nara2009,Ni2012,NishWaka2014,SunWang2007}) are quite different: roughly speaking, this is due to the strong influence of the classical damping terms on the long time behavior of the solutions to the linearized model and, in turn, on the range of the powers to which is associated the global in-time prolongability of the small-data solutions to the semilinear problem. As a consequence of this different expression for the critical curves, which is the counterpart for weakly coupled systems of what happens with the Strauss and Fujita exponents for the corresponding single semilinear equations with power-nonlinearities, it arises quite naturally the question on the critical curve for a system in which a wave and a damped wave equations are coupled through the semilinear terms. The Cauchy problem for this kind of system is named Nakao problem, after the author of \cite{N16,N18}, and it has been recently investigated in several directions: from the first results for the Nakao problem with power nonlinearities \cite{Chen-Reissig=2021, Kita-Kusaba=2022, Wakasugi=2017}, to the Nakao problem with derivative type nonlinearities \cite{Chen=2022, Liu=2026, Palmieri-Takamura=2023}, and lately to Nakao-type problems with time-dependent coefficients for the damping terms \cite{Li-Palmieri=2025, Li-Palmieri2=2025}.

In the present paper, we consider in some sense a Nakao-type problem when the space variable runs in a compact Lie group $\mathbb{G}$. We stress that in our model in \eqref{Eq-Main-Systems} we can consider both the Klein-Gordon equation with dominant damping, i.e. for $b^2\geqslant 4m^2 >0$ (hereafter the \emph{damped Klein-Gordon case}) and the classical damped equation, i.e. for $b^2>4m^2=0$ (hereafter the \emph{massless case}), obtaining the same blow-up result but quite different lifespan estimates.

We point out that, while for the Euclidean models we have a critical exponent for the corresponding semilinear equation that is the threshold value for the exponent in the nonlinearity that separates the blow-up range from the small-data solutions' global existence range, in the case of a compact Lie group the situation is rather bifurcated: in the literature (cf. \cite{Bhardwaj-Kumar-Mondal=2024, Dasgupta-Kumar-Shyam-2023, Palmieri-JDE=2021, Palmieri-Pseudo=2021, Palmieri-JFAA=2022, Ruzhansky-Yess=2024}) there are examples of Cauchy problems on compact Lie groups for which a blow-up result holds for any $p>1$ but also examples where a global existence result for small-data solutions holds for any $p>1$. This dichotomous scenario for semilinear evolution models on a compact Lie group $\mathbb{G}$ is essentially due to the fact that the global dimension of the group is 0 (being $\mathbb{G}$ bounded). For further details on the global dimension of a Lie group, we refer to \cite{Var1988} or \cite[Section II.4]{DER2003}. 

For the model \eqref{Eq-Main-Systems} it will be interesting to see that the solutions blow up in finite time for any $p,q>1$ (under suitable integral sign assumptions for the Cauchy data), in spite of the fact that for the single semilinear damped Klein-Gordon equation we have the global existence of small-data solutions for any power greater than 1 (see \cite{Palmieri-JFAA=2022}).

Our main theorems will consist of two kinds of results: on the one hand, blow-up results obtained by using the approach for semilinear wave equations, introduced in \cite{K1980} and later developed (among the others) in \cite{G1981, Sc1985, Si1984, Ta2015, Zh2006}, that provides as byproducts the upper bound estimates for the lifespan as well; on the other hand, local in-time existence results on $L^2(\mathbb{G})$-based Sobolev spaces that yield lower bound estimates for the lifespan.

In the method employed for the blow-up results, we apply the so-called \emph{slicing procedure}, a technique introduced in \cite{AgKuTa2000}, in the formulation from \cite{Chen-Palmieri=2020, Palmieri-Takamura=2022, Palmieri-Takamura=2023}. For the existence results, we employ a standard contraction argument based on the $L^2(\mathbb{G})$ estimates for the corresponding linearized problem with a suitable loss of decay in the damped Klein-Gordon case. These estimates are obtained by using the Fourier series on compact Lie groups. The core of this theory is the Peter-Weyl theorem \cite{PW27}, that allows to explicitly describe the unitary dual $\widehat{\mathbb{G}}$ and to derive a Plancherel formula for the Fourier series on $\mathbb{G}$. For a modern description of this topic we refer to \cite{FisherRuzhansky2016, RuzhanskyTurunen2010}.
Finally, we underline that, the employment of the Fourier series for studying the well-posedness of linear wave equations on a compact Lie group was introduced in \cite{GarettoRuzhansky2015}.

\paragraph{\large Notation:}The positive constants $c$ and $C$ may change from line to line. We write $f\lesssim g$ if there exists a positive constant $C$ such that $f\leqslant Cg$. $L^{r}(\G)$ denotes the space of $r$-summable functions on $\G$ with respect to the normalized Haar measure for $r\geqslant 1$ 
(respectively, essentially bounded functions for $r=\infty$). For $\sigma>0$ and $r> 1$ 
the Sobolev space $H^{\sigma,r}_{\ml{L}}(\G)$ is defined as the space
\begin{align*}
H^{\sigma,r}_{\ml{L}}(\G):=\left\{f\in L^{r}(\G):\ (-\ml{L})^{\sigma/2}f\in L^{r}(\G) \right\}
\end{align*}
equipped with the norm $\|f\|_{H^{\sigma,r}_{\ml{L}}(\G)}:=\|f\|_{L^{r}(\G)}+\|(-\ml{L})^{\sigma/2}f\|_{L^{r}(\G)}$. As usual, the Hilbert space $H^{\sigma,2}_{\ml{L}}(\G)$ is simply denoted by $H^{\sigma}_{\ml{L}}(\G)$.

\section{Main results}\setcounter{equation}{0}\label{Section-Main-Result}
\hspace{5mm}Before stating our blow-up results, let us introduce the notion of energy solutions for the weakly coupled system \eqref{Eq-Main-Systems} 
that we are going to employ in our theorems.
\begin{defn}\label{Defn-Energy-Solutions} Let $b>0$ and $m^2\geqslant 0$.
Let $(u_0,u_1),(v_0,v_1)\in H^1_{\ml{L}}(\G)\times L^2(\G)$. We say that $(u,v)$ is an \emph{energy solution} on $[0,T)$ to the weakly coupled system \eqref{Eq-Main-Systems} if
\begin{align*}
u&\in\ml{C}\big([0,T),H^1_{\ml{L}}(\G)\big)\cap\ml{C}^1\big([0,T),L^2(\G)\big)\cap L^q_{\mathrm{loc}}\big([0,T)\times \G\big),\\
v&\in\ml{C}\big([0,T),H^1_{\ml{L}}(\G)\big)\cap\ml{C}^1\big([0,T),L^2(\G)\big)\cap L^p_{\mathrm{loc}}\big([0,T)\times \G\big),
\end{align*}
and if $(u,v)$ fulfills the following integral relations:
\begin{align*}
\int_0^t\int_{\G}|v(s,x)|^p\psi(s,x)\,\mathrm{d}x\,\mathrm{d}s&=\int_{\G}\partial_tu(t,x)\psi(t,x)\,\mathrm{d}x-\int_{\G} u(t,x)\partial_s\psi(t,x)\,\mathrm{d}x+b\int_{\G}u(t,x)\psi(t,x)\,\mathrm{d}x\\
&\quad-\varepsilon\int_{\G}u_1(x)\psi(0,x)\,\mathrm{d}x+\varepsilon\int_{\G}u_0(x)\partial_s\psi(0,x)\,\mathrm{d}x-\varepsilon b\int_{\G}u_0(x)\psi(0,x)\,\mathrm{d}x\\
&\quad+\int_0^t\int_{\G}u(s,x)\big(\partial_s^2\psi(s,x)-\ml{L}\psi(s,x)-b\,\partial_s\psi(s,x)+m^2\psi(s,x)\big)\,\mathrm{d}x\,\mathrm{d}s
\end{align*}
and
\begin{align*}
\int_0^t\int_{\G}|u(s,x)|^q\varphi(s,x)\,\mathrm{d}x\,\mathrm{d}s&=\int_{\G}\partial_tv(t,x)\varphi(t,x)\,\mathrm{d}x-\int_{\G}v(t,x)\partial_s\varphi(t,x)\,\mathrm{d}x-\varepsilon\int_{\G}v_1(x)\varphi(0,x)\,\mathrm{d}x\\
&\quad+\varepsilon\int_{\G}u_0(x)\partial_s\varphi(0,x)\,\mathrm{d}x+\int_0^t\int_{\G}v(s,x)\big(\partial_s^2\varphi(s,x)-\ml{L}\varphi(s,x)\big)\,\mathrm{d}x\,\mathrm{d}s 
\end{align*}
for any $\psi,\varphi\in\ml{C}_0^{\infty}([0,T)\times\G)$ and any $t\in(0,T)$.
\end{defn}

\begin{remark} Let $\varepsilon>0$. If $(u,v)$ is an energy solution to \eqref{Eq-Main-Systems} according to Definition \ref{Defn-Energy-Solutions}, then the lifespan of $(u,v)$ is defined by the quantity
\begin{align*}
T(\varepsilon) := \sup \big\{\tau >0: (u,v) \mbox{ is an energy solution to \eqref{Eq-Main-Systems} on } [0,\tau) \big\}.
\end{align*} If $T(\varepsilon)=+\infty$ we call $(u,v)$ a global in-time energy solution to \eqref{Eq-Main-Systems}. 

We stress that, given $\varepsilon>0$, the solution $(u,v)$ depends on $\varepsilon$, however, for the sake of simplicity, we shall avoid the notation $(u_\varepsilon,v_\varepsilon)$, keeping the dependence of the solution on $\varepsilon$ implicit. 
\end{remark}

\begin{theorem}\label{Thm-Blow-up}
Let $\G$ be a compact Lie group. Let $\varepsilon>0$ and $b^2\geqslant 4m^2>0$. Let $u_0,v_0 \in H^1_{\ml{L}}(\G)$, $u_1,v_1 \in L^2(\G)$ be non-negative functions such that each pair $(u_0,u_1), (v_0,v_1)$ has at least one nontrivial component. 

Let $p,q>1$ and let $(u,v)$ be an energy solution to the weakly coupled system \eqref{Eq-Main-Systems} according to Definition \ref{Defn-Energy-Solutions} with lifespan $T(\varepsilon)$. Then, there exists a positive constant $\varepsilon_0=\varepsilon_0(u_0,u_1,v_0,v_1,b,m^2,p,q)>0$ such that for any $\varepsilon\in(0,\varepsilon_0]$ the energy solution $(u,v)$ blows up in finite time. 

Furthermore, $T(\varepsilon)$ satisfies the following upper bound estimates:
\begin{align}\label{Est-2}
		T(\varepsilon)\leqslant \begin{cases}
			C\varepsilon^{-\frac{pq-1}{pq+1}}&\mbox{if}\ \ v_1\neq0,\\
			C\varepsilon^{-\frac{pq-1}{2}}&\mbox{if}\ \ v_1=0,
		\end{cases}
\end{align}
where the constant $C>0$ is independent of $\varepsilon$.
\end{theorem}

\begin{theorem}\label{Thm-Blow-up massless}
Let $\G$ be a compact Lie group. Let $\varepsilon>0$ and $b^2>4m^2=0$. Let $u_0,v_0 \in H^1_{\ml{L}}(\G)$, $u_1,v_1 \in L^2(\G)$ be non-negative functions such that each pair $(u_0,u_1), (v_0,v_1)$ has at least one nontrivial component. 

Let $p,q>1$ and let $(u,v)$ be an energy solution to the weakly coupled system \eqref{Eq-Main-Systems} according to Definition \ref{Defn-Energy-Solutions} with lifespan $T(\varepsilon)$. Then, there exists a positive constant $\varepsilon_0=\varepsilon_0(u_0,u_1,v_0,v_1,b,m^2,p,q)>0$ such that for any $\varepsilon\in(0,\varepsilon_0]$ the energy solution $(u,v)$ blows up in finite time. 

Furthermore, $T(\varepsilon)$ satisfies the following upper bound estimates:
\begin{align*}
		T(\varepsilon)\leqslant \begin{cases}
			C\varepsilon^{-\frac{pq-1}{\max\{2p+1,pq+q+1\}}}&\mbox{if}\ \ v_1\neq0,\\
			C\varepsilon^{-\frac{pq-1}{\max\{2p+1,q+2\}}}&\mbox{if}\ \ v_1=0,
		\end{cases}
\end{align*}
where the constant $C>0$ is independent of $\varepsilon$.
\end{theorem}

After the blow-up results, we focus now on the existence of the local in-time mild solutions and the corresponding lower bound estimates for the lifespan.

\begin{theorem}\label{Thm-Local}
	Let $\G$ be a compact, connected Lie group. Let us assume that the topological dimension $n$ of $\G$ satisfies $n\geqslant 3$. Let  $b^2\geqslant 4m^2>0$, $\varepsilon>0$ and $(u_0,u_1),(v_0,v_1)\in H^1_{\ml{L}}(\G)\times L^2(\G)$. Let $p,q>1$ such that $p,q\leqslant\frac{n}{n-2}$. 
	
If $T>0$ satisfies
\begin{align*}
T\lesssim \begin{cases}
\varepsilon^{-\frac{pq-1}{pq+1}}&\mbox{if}\ \ v_1\neq0,\\
\varepsilon^{-\frac{pq-1}{2}}&\mbox{if}\ \ v_1=0 \ \  \mbox{and}\ \ p<2-\frac{1}{q}, 
\end{cases}
\end{align*}
then the weakly coupled system \eqref{Eq-Main-Systems} admits a uniquely determined mild solution
\begin{align*}
(u,v)\in \left(\ml{C}\big([0,T],H^1_{\ml{L}}(\G)\big)\cap \ml{C}^1\big([0,T],L^2(\G)\big)\right)^2.
\end{align*}
In particular, the lifespan $T(\varepsilon)$ satisfies the following lower bound estimates:
\begin{align}\label{Est-01}
T(\varepsilon)\geqslant \begin{cases}
c\,\varepsilon^{-\frac{pq-1}{pq+1}}&\mbox{if}\ \ v_1\neq0,\\
c\,\varepsilon^{-\frac{pq-1}{2}}&\mbox{if}\ \ v_1=0\ \  \mbox{and}\ \  p<2-\frac{1}{q},
\end{cases}
\end{align}
where the constant $c>0$ is independent of $\varepsilon$.
\end{theorem}

\begin{remark}
 In the statement of Theorem \ref{Thm-Local}, the assumptions on the connectedness of the group $\G$, on the topological dimension $n$ ($n\geqslant 3$) and on the powers of the nonlinear terms (the upper bounds $p,q\leqslant\frac{n}{n-2}$ ) are made in order to apply the Gagliardo-Nirenberg type inequality proved in \cite[Remark 1.7]{Ruzhansky-Yess=2024}.
\end{remark}

\begin{remark}
Combining \eqref{Est-2} and \eqref{Est-01}, we obtain the sharp lifespan estimates
\begin{align*}
\begin{cases}
	c\,\varepsilon^{-\frac{pq-1}{pq+1}}\leqslant T(\varepsilon)\leqslant C\varepsilon^{-\frac{pq-1}{pq+1}}&\mbox{if}\ \ v_1\neq0,\\
	c\,\varepsilon^{-\frac{pq-1}{2}}\leqslant T(\varepsilon)\leqslant C\varepsilon^{-\frac{pq-1}{2}}&\mbox{if}\ \ v_1=0\ \  \mbox{and}\ \ p<2-\frac{1}{q},
\end{cases}
\end{align*}
for local in-time solutions to the weakly coupled system \eqref{Eq-Main-Systems} with $b^2\geqslant 4m^2 >0$. Therefore, the non-triviality of $v_1$ plays a crucial role in determining the lifespan of the local solution $(u,v)$.
\end{remark}

\begin{theorem}\label{Thm-Local massless case}
	Let $\G$ be a compact, connected Lie group. Let us assume that the topological dimension $n$ of $\G$ satisfies $n\geqslant 3$.  Let $b^2>4m^2=0$, $\varepsilon>0$ and $(u_0,u_1),(v_0,v_1)\in H^1_{\ml{L}}(\G)\times L^2(\G)$. Let $p,q>1$ such that $p,q\leqslant\frac{n}{n-2}$. 
	
If $T>0$ satisfies
\begin{align*}
T \lesssim \begin{cases}
\varepsilon^{-\min\left\{\frac{p-1}{p+1},q-1 \right\}}&\mbox{if}\ \ v_1\neq0,\\
\varepsilon^{-\min\left\{p-1,\frac{q-1}{2} \right\}}&\mbox{if}\ \ v_1=0,
\end{cases}
\end{align*}
then, the weakly coupled system \eqref{Eq-Main-Systems} admits a uniquely determined mild solution
\begin{align*}
(u,v)\in\ml{C}\big([0,T],H^1_{\ml{L}}(\G)\big)\times \left(\ml{C}\big([0,T],H^1_{\ml{L}}(\G)\big)\cap \ml{C}^1\big([0,T],L^2(\G)\big)\right).
\end{align*}
In particular, the lifespan $T(\varepsilon)$ satisfies the following lower bound estimates:
\begin{align*}
T(\varepsilon)\geqslant \begin{cases}
c\,\varepsilon^{-\min\left\{\frac{p-1}{p+1},q-1 \right\}}&\mbox{if}\ \ v_1\neq0,\\
c\,\varepsilon^{-\min\left\{p-1,\frac{q-1}{2} \right\}}&\mbox{if}\ \ v_1=0,
\end{cases}
\end{align*}
where the constant $c>0$ is independent of $\varepsilon$.
\end{theorem}

\section{Proof of Theorem \ref{Thm-Blow-up}} 
\label{Section-Blow-up-damping-mass}\setcounter{equation}{0}
\hspace{5mm}Let $(u,v)$ be a local in-time energy solution to the weakly coupled system \eqref{Eq-Main-Systems} with $m^2>0$ according to Definition \ref{Defn-Energy-Solutions} with the lifespan $T$. Let us fix $t\in(0,T)$. We can choose two bump functions $\psi,\varphi\in\ml{C}_0^{\infty}([0,T)\times\G)$ such that $\psi\equiv 1\equiv \varphi$ on $[0,t]\times\G$. Therefore,
\begin{align}\label{Eq-from-energy-1}
\int_0^t\int_{\G}|v(s,x)|^p\,\mathrm{d}x\,\mathrm{d}s&=\int_{\G}\partial_tu(t,x)\,\mathrm{d}x+b\int_{\G}u(t,x)\,\mathrm{d}x+m^2\int_0^t\int_{\G}u(s,x)\,\mathrm{d}x\,\mathrm{d}s\notag\\
&\quad-\varepsilon\int_{\G}\big(u_1(x)+bu_0(x)\big)\,\mathrm{d}x
\end{align}
and
\begin{align}\label{Eq-from-energy-2}
\int_0^t\int_{\G}|u(s,x)|^q\,\mathrm{d}x\,\mathrm{d}s=\int_{\G}\partial_tv(t,x)\,\mathrm{d}x-\varepsilon\int_{\G}v_1(x)\,\mathrm{d}x.
\end{align}

\subsection{Iteration frames and first lower bounds}
\hspace{5mm}Let us introduce two time-dependent functionals
\begin{align}\label{Functionals}
	\ml{U}(t):=\int_{\G}u(t,x)\,\mathrm{d}x\ \ \mbox{and}\ \ \ml{V}(t):=\int_{\G}v(t,x)\,\mathrm{d}x.
\end{align} Differentiating \eqref{Eq-from-energy-1} with respect to $t$, it yields
\begin{align*}
\ml{U}''(t)+b\ml{U}'(t)+m^2\ml{U}(t)=\int_{\G}|v(t,x)|^p\,\mathrm{d}x\geqslant|\ml{V}(t)|^p,
\end{align*}
where we used Jensen's inequality and the normalization for the Haar measure on $\mb{G}$. The differential operator on the left-hand side of the previous equation can be factorized as
\begin{align*}
\mathrm{e}^{-k_2t}\frac{\mathrm{d}}{\mathrm{d}t}\Big(\mathrm{e}^{(k_2-k_1)t}\frac{\mathrm{d}}{\mathrm{d}t}\big(\mathrm{e}^{k_1t}\ml{U}(t)\big)\Big)=\ml{U}''(t)+(k_1+k_2)\ml{U}'(t)+k_1k_2\ml{U}(t)
\end{align*}
thanks to the condition $b^2\geqslant 4m^2$, where $k_1,k_2$ solve the system
\begin{align*}
\begin{cases}
k_1+k_2=b, \\
k_1k_2=m^2,
\end{cases}
\end{align*} namely, we have
\begin{align*}
k_{1/2}:=\frac{1}{2}\big(b\pm\sqrt{b^2-4m^2}\,\big)>0.
\end{align*}
 Thus, we arrive at
\begin{align}\label{Ineq-01}
	\mathrm{e}^{-k_2t}\frac{\mathrm{d}}{\mathrm{d}t}\Big(\mathrm{e}^{(k_2-k_1)t}\frac{\mathrm{d}}{\mathrm{d}t}\big(\mathrm{e}^{k_1t}\ml{U}(t)\big)\Big)\geqslant|\ml{V}(t)|^p.
\end{align}
\begin{remark}
In the case with dominant damping $b^2>4m^2$, there are two distinct real roots $k_{1/2}$ to the quadratic equation $k^2-bk+m^2=0$. On the contrary, in the case with dominant mass $b^2<4m^2$, the complex conjugate roots to the quadratic equation lead to oscillations. Finally, in the balanced case $b^2=4m^2$, there is a double root $k_1=k_2=\frac{b}{2}$. 
\end{remark}

Next, we are going to derive an iteration frame for the functional $\ml{V}(t)$. Integrating \eqref{Ineq-01} twice implies
\begin{align}
\ml{U}(t) & \geqslant \mathrm{e}^{-k_1t}\ml{U}(0)+\mathrm{e}^{-k_1t}\int_0^t\mathrm{e}^{(k_1-k_2)s}\,\mathrm{d}s \ \big(k_1\ml{U}(0)+\ml{U}'(0)\big) \notag \\& \quad +\mathrm{e}^{-k_1t}\int_0^t\mathrm{e}^{(k_1-k_2)s}\int_0^s\mathrm{e}^{k_2\tau}|\ml{V}(\tau)|^p\,\mathrm{d}\tau\,\mathrm{d}s, \label{ast}
\end{align}
where
\begin{align*}
\mathrm{e}^{-k_1t}\int_0^t\mathrm{e}^{(k_1-k_2)s}\,\mathrm{d}s=\begin{cases}
\displaystyle{\frac{\mathrm{e}^{-k_2t}-\mathrm{e}^{-k_1t}}{k_1-k_2}}&\mbox{if}\ \ b^2>4m^2,\\[0.8em]
t\,\mathrm{e}^{-\frac{b}{2}t}&\mbox{if}\ \ b^2=4m^2.
\end{cases}
\end{align*}
Let us introduce the function
\begin{align}\label{Lower-Bound-U-1}
\ml{U}^{\mathrm{lin}}(t) := \begin{cases}
\displaystyle{\frac{k_1\mathrm{e}^{-k_2t}-k_2\mathrm{e}^{-k_1t}}{k_1-k_2}\ml{U}(0)+\frac{\mathrm{e}^{-k_2t}-\mathrm{e}^{-k_1t}}{k_1-k_2}\ml{U}'(0)}&\mbox{if}\ \ b^2>4m^2,\\[1em]
\mathrm{e}^{-\frac{b}{2}t}\big((1+\tfrac{b}{2}t)\ml{U}(0)+t\ml{U}'(0)\big)&\mbox{if}\ \ b^2=4m^2.
\end{cases}
\end{align}  Since we can rewrite \eqref{ast} as
\begin{align}\label{ast 2}
\ml{U}(t) & \geqslant \ml{U}^{\mathrm{lin}}(t) +\mathrm{e}^{-k_1t}\int_0^t\mathrm{e}^{(k_1-k_2)s}\int_0^s\mathrm{e}^{k_2\tau}|\ml{V}(\tau)|^p\,\mathrm{d}\tau\,\mathrm{d}s,
\end{align} and the non-negativity of the initial data $u_0,u_1$ guarantees that the function $\ml{U}^{\mathrm{lin}}$ is non-negative, we obtained the nonlinear integral inequality
\begin{align}\label{Frame-U}
\ml{U}(t)\geqslant \mathrm{e}^{-k_1t}\int_0^t\mathrm{e}^{(k_1-k_2)s}\int_0^s\mathrm{e}^{k_2\tau}|\ml{V}(\tau)|^p\,\mathrm{d}\tau\,\mathrm{d}s.
\end{align}
On the other hand, from \eqref{Eq-from-energy-2} we have
\begin{align*}
\ml{V}(t)\geqslant\ml{V}(0)+\ml{V}'(0)\,t+\int_0^t\int_0^s|\ml{U}(\tau)|^q\,\mathrm{d}\tau\,\mathrm{d}s.
\end{align*}
Thanks to our assumption on  $v_0,v_1$, one derives a first lower bound for $\ml{V}$
\begin{align}\label{Lower-Bound-V-1}
\ml{V}(t)\geqslant C_r\varepsilon\,t^r\ \ \mbox{with}\ \ C_r:=\int_{\G}v_r(x)\,\mathrm{d}x, \ 
\end{align} where either $r=0$ if $v_1\equiv 0$ or $r=1$ otherwise, and the following  nonlinear integral inequality:
\begin{align}\label{Frame-V}
\ml{V}(t)\geqslant \int_0^t\int_0^s|\ml{U}(\tau)|^q\,\mathrm{d}\tau\,\mathrm{d}s.
\end{align}
The inequalities in \eqref{Frame-U}, \eqref{Frame-V} provide the iteration frame for $\ml{V}$.

\subsection{Iteration argument}
\hspace{5mm}From \eqref{ast 2} we see that $\ml{U}(t)\geqslant \ml{U}^{\mathrm{lin}}(t)$.
Due to the exponentially decay of $\ml{U}^{\mathrm{lin}}(t)$ in \eqref{Lower-Bound-U-1}, the time-dependent functional that we consider to prove the blow-up result is the space average $\ml{V}(t)$, which has on the contrary a polynomially growing or bounded from below lower bound thanks to \eqref{Lower-Bound-V-1}.

To overcome some difficulties from the unbounded exponential multipliers $\mathrm{e}^{(k_1-k_2)s}$ and $\mathrm{e}^{k_2\tau}$ in the $s$-integral and in the $\tau$-integral in \eqref{Frame-U}, respectively, we are going to apply a 2-step slicing procedure adapting the main ideas from \cite{Chen-Palmieri=2020,Chen-Reissig=2021,Palmieri-Takamura=2022} to our problem. More specifically, we introduce the sequence $\{L_j\}_{j\in\mb{N}_0}$ that characterizes the slicing procedure by defining
\begin{align}\label{Defn-Lj}
L_j:=\prod\limits_{k=0}^{j}\ell_k\ \ \mbox{for any}\ \ j\in\mb{N}_0,
\end{align}
where the elements of the sequence $\{\ell_k\}_{k\in\mb{N}_0}$ are
\begin{align*}
\ell_0:=\frac{1}{k_2}\ \ \mbox{and}\ \  \ell_k:=1+(pq)^{-\frac{k}{2}}\ \ \mbox{for any}\ \ k\geqslant 1.
\end{align*}
Notice that $\ell_k>1$ for any $k\geqslant 1$, hence, the sequence $\{L_j\}_{j\in\mb{N}_0}$ is strictly increasing. Thanks to the choice of $\{\ell_k\}_{k\geqslant 1}$, the series $\sum_{k=1}^{+\infty}\ln\ell_k$ is convergent, and this is equivalent to the convergence of the following infinite product:
\begin{align*}
L:=\prod\limits_{j=0}^{+\infty}\ell_j\in\mb{R}_+.
\end{align*}

We will determine a sequence of lower bounds for $\ml{V}(t)$ via \eqref{Frame-U} and \eqref{Frame-V} in an iterative way. To be specific, we show that
\begin{align}\label{Seq-V}
\ml{V}(t)\geqslant D_j(t-L_{2j})^{\alpha_j}\ \ \mbox{for any}\ \ t\geqslant L_{2j},
\end{align}
where $\{D_j\}_{j\in\mb{N}_0}$ and $\{\alpha_j\}_{j\in\mb{N}_0}$ are sequences of non-negative real numbers to be determined iteratively later. From \eqref{Lower-Bound-V-1} we have the validity of base case $j=0$ with $D_0:=C_r\varepsilon$ and $\alpha_0:=r$.

Next, we prove the inductive step. Assuming that \eqref{Seq-V} holds for some $j\geqslant 0$, we will prove it for $j+1$. Plugging \eqref{Seq-V} into \eqref{Frame-U}, one has
\begin{align*}
\ml{U}(t)\geqslant D_j^p\,\mathrm{e}^{-k_1t}\int_{L_{2j}}^t\mathrm{e}^{(k_1-k_2)s}\int_{L_{2j}}^s\mathrm{e}^{k_2\tau}(\tau-L_{2j})^{\alpha_jp}\,\mathrm{d}\tau\,\mathrm{d}s.
\end{align*}
For $t\geqslant L_{2j+1}$ and $s\in [ L_{2j+1},t]$, we can shrink the domain of integration from $[L_{2j},s]$ to $[\frac{s}{\ell_{2j+1}},s]$, that is,
\begin{align*}
\ml{U}(t)&\geqslant D_j^p\,\mathrm{e}^{-k_1t}\int_{L_{2j+1}}^t\mathrm{e}^{(k_1-k_2)s}\int_{\tfrac{s}{\ell_{2j+1}}}^s\mathrm{e}^{k_2\tau}(\tau-L_{2j})^{\alpha_jp}\,\mathrm{d}\tau\,\mathrm{d}s\\
&\geqslant\frac{D_j^p}{k_2\ell_{2j+1}^{\alpha_jp}}\,\mathrm{e}^{-k_1t}\int_{L_{2j+1}}^t\mathrm{e}^{k_1s}(s-L_{2j+1})^{\alpha_jp}\left(1-\mathrm{e}^{-k_2(1-\frac{1}{\ell_{2j+1}})s}\right)\mathrm{d}s \\
&\geqslant\frac{D_j^p}{k_2\ell_{2j+1}^{\alpha_jp}} \left(1-\mathrm{e}^{-(\ell_{2j+1}-1)}\right) \mathrm{e}^{-k_1t}\int_{L_{2j+1}}^t\mathrm{e}^{k_1s}(s-L_{2j+1})^{\alpha_jp}\,\mathrm{d}s,
\end{align*}
where we used $k_2L_{2j}\geqslant k_2\ell_0=1$ and the monotonicity of the exponential function in the last inequality.
Then, using the inequality $1-\mathrm{e}^{-y}\geqslant y-\frac{y^2}{2}$ for any $y\geqslant0$, it results
\begin{align}
1-\mathrm{e}^{-(\ell_{2j+1}-1)}&\geqslant (\ell_{2j+1}-1)\big(1-\tfrac{1}{2}(\ell_{2j+1}-1)\big) \notag \\
&\geqslant (pq)^{-(2j+1)}\big((pq)^{\frac{2j+1}{2}}-\tfrac{1}{2}\big) \notag \\
&\geqslant(pq)^{-(2j+1)}\big(pq-\tfrac{1}{2}\big). \label{ast 3}
\end{align}
Similarly, for $t\geqslant L_{2j+2}$, by shrinking the $s$-domain of integration from $[L_{2j+1},t]$ to $[\frac{t}{\ell_{2j+2}},t]$, we find
\begin{align}\label{Bound-1}
	\ml{U}(t)
	&\geqslant\frac{D_j^p(pq-\tfrac{1}{2})}{k_2\ell_{2j+1}^{\alpha_jp}(pq)^{2j+1}}\,\mathrm{e}^{-k_1t}\int_{\tfrac{t}{\ell_{2j+2}}}^t\mathrm{e}^{k_1s}(s-L_{2j+1})^{\alpha_jp}\,\mathrm{d}s\notag\\
	&\geqslant \frac{D_j^p(pq-\tfrac{1}{2})}{k_1k_2(\ell_{2j+1}\ell_{2j+2})^{\alpha_jp}(pq)^{2j+1}}(t-L_{2j+2})^{\alpha_jp}\Big(1-\mathrm{e}^{-k_1(1-\frac{1}{\ell_{2j+2}})t}\Big)\notag\\
	&\geqslant \frac{D_j^p(pq-\tfrac{1}{2})^2}{k_1k_2(\ell_{2j+1}\ell_{2j+2})^{\alpha_jp}(pq)^{4j+3}}(t-L_{2j+2})^{\alpha_jp},
\end{align}
in which we estimated
\begin{align*}
1-\mathrm{e}^{-k_1(1-\frac{1}{\ell_{2j+2}})t} & =1-\mathrm{e}^{-k_1(\ell_{2j+2}-1)\frac{t}{\ell_{2j+2}}} \\
& \geqslant 1-\mathrm{e}^{-k_1(\ell_{2j+2}-1)L_{2j+1}} \\ & \geqslant 1-\mathrm{e}^{-(\ell_{2j+2}-1)}
\end{align*}
by using $k_1L_{2j+1}\geqslant k_1\ell_0\geqslant 1$ and \eqref{ast 3}. Finally, for $t\geqslant L_{2j+2}$, we plug the lower bound \eqref{Bound-1} into \eqref{Frame-V}, arriving at
\begin{align*}
\ml{V}(t)&\geqslant \frac{D_j^{pq}(pq-\tfrac{1}{2})^{2q}}{(k_1k_2)^q(\ell_{2j+1}\ell_{2j+2})^{\alpha_jpq}(pq)^{(4j+3)q}}\int_{L_{2j+2}}^t\int_{L_{2j+2}}^s(\tau-L_{2j+2})^{\alpha_jpq}\,\mathrm{d}\tau\,\mathrm{d}s\\
&\geqslant \frac{D_j^{pq}(pq-\tfrac{1}{2})^{2q}}{(k_1k_2)^q(\ell_{2j+1}\ell_{2j+2})^{\alpha_jpq}(\alpha_jpq+1)(\alpha_jpq+2)(pq)^{(4j+3)q}}(t-L_{2j+2})^{\alpha_jpq+2},
\end{align*}
which is exactly \eqref{Seq-V} for $j+1$ provided that
\begin{align*}
\alpha_{j+1}:=\alpha_jpq+2\ \ \mbox{and}\ \ D_{j+1}:=\frac{(pq-\tfrac{1}{2})^{2q}}{(k_1k_2)^q(\ell_{2j+1}\ell_{2j+2})^{\alpha_jpq}(\alpha_jpq+1)(\alpha_jpq+2)(pq)^{(4j+3)q}}D_j^{pq}.
\end{align*}
\subsection{Upper bound estimates for the lifespan} \label{Subsection 3.3}
\hspace{5mm}In the previous subsection, we determined a sequence of lower bound estimates for the functional $\ml{V}(t)$. Now, we are going to show that the $j$-dependent lower bound blows up as $j\to+\infty$ when $t$ is greater than a certain $\varepsilon$-dependent threshold. This will demonstrate the desired blow-up result and, as a byproduct, will provide the upper bound estimates for the lifespan.

The recursive relation $\alpha_j=\alpha_{j-1}pq+2$ yields
\begin{align}\label{alphaj}
\alpha_j=2\big(1+pq+\cdots+(pq)^{j-1}\big)+\alpha_0(pq)^j=\big(\tfrac{2}{pq-1}+r\big)(pq)^j-\tfrac{2}{pq-1},
\end{align}
and, therefore,
\begin{align*}
(\alpha_jpq+1)(\alpha_jpq+2)\leqslant \alpha_{j+1}^2\leqslant \big(\tfrac{2}{pq-1}+r\big)^2(pq)^{2j+2}=:M_0(pq)^{2j+2}.
\end{align*}
Since
\begin{align*}
\lim\limits_{j\to+\infty}(\ell_{2j+1}\ell_{2j+2})^{\alpha_{j+1}}&=\lim\limits_{j\to+\infty}\exp\big(\alpha_{j+1}(\ln\ell_{2j+1}+\ln\ell_{2j+2})\big)\\
&=\lim\limits_{j\to+\infty}\exp\left\{\big(\tfrac{2}{pq-1}+r\big)(pq)^{j+1}\left[\ln\big(1+(pq)^{-(j+\frac{1}{2})}\big)+\ln\big(1+(pq)^{-(j+1)}\big)\right]\right\}\\
&=\exp\big[\big(\tfrac{2}{pq-1}+r\big)(1+\sqrt{pq}\,)\big],
\end{align*}
there exists a uniform (i.e. independent of $j$) constant $M_1=M_1(p,q,r)>0$ satisfying $$(\ell_{2j+1}\ell_{2j+2})^{\alpha_{j+1}}\leqslant M_1$$ for any $j\in\mb{N}_0$. Summarizing, we may estimate
\begin{align*}
D_{j+1}\geqslant\frac{(pq-\frac{1}{2})^{2q}(pq)^q}{(k_1k_2)^qM_0M_1}(pq)^{-(4q+2)(j+1)}D_j^{pq}=:M_2(pq)^{-(4q+2)(j+1)}D_j^{pq}
\end{align*}
for a suitable constant $M_2=M_2(b,m^2,p,q,v_0,v_1)>0$. Then, applying the logarithmic function to both sides of the previous inequality and using iteratively the resultant inequality, we deduce
\begin{align*}
\ln D_j&\geqslant(pq)\ln D_{j-1}-(4q+2)j\ln(pq)+\ln M_2\\
&\geqslant\cdots\geqslant (pq)^j\ln D_0-(4q+2)\left(\ \sum\limits_{k=0}^{j-1}(j-k)(pq)^k\right)\ln(pq)+\left(\ \sum\limits_{k=0}^{j-1}(pq)^k\right)\ln M_2.
\end{align*}
According to the identity
\begin{align*}
\sum\limits_{k=0}^{j-1}(j-k)(pq)^k=\frac{1}{pq-1}\left(\frac{(pq)^{j+1}-pq}{pq-1}-j\right),
\end{align*}
one concludes
\begin{align}\label{Dj}
\ln D_j&\geqslant(pq)^j\left(\ln D_0-\frac{(4q+2)pq\ln(pq)}{(pq-1)^2}+\frac{\ln M_2}{pq-1}\right)+\frac{(4q+2)pq\ln(pq)}{(pq-1)^2}+\frac{(4q+2)\ln(pq)}{pq-1}j-\frac{\ln M_2}{pq-1}\notag\\
&\geqslant(pq)^j\left(\ln D_0-\frac{(4q+2)pq\ln(pq)}{(pq-1)^2}+\frac{\ln M_2}{pq-1}\right)
\end{align}
for any $j\geqslant j_0$, where $j_0=j_0(b,m^2,p,q,v_0,v_1)\in\mb{N}_0$ is the smallest integer such that
\begin{align*}
j_0\geqslant\frac{\ln M_2}{(4q+2)\ln(pq)}-\frac{pq}{pq-1}.
\end{align*}
Recalling that $L_{2j}\uparrow L$ as $j\to+\infty$, and combining \eqref{alphaj}, \eqref{Dj} in \eqref{Seq-V}, the following lower bound estimate holds for any $t\geqslant L$ and any $j\geqslant j_0$: 
\begin{align*}
\ml{V}(t)
&\geqslant\exp\left[(pq)^j\left(\ln D_0-\frac{(4q+2)pq\ln(pq)}{(pq-1)^2}+\frac{\ln M_2}{pq-1}+\left(\frac{2}{pq-1}+r\right)\ln(t-L)\right) \right](t-L)^{-\frac{2}{pq-1}}.
\end{align*}
Choosing $t\geqslant 2L$ so that $t-L\geqslant\frac{t}{2}$, and denoting
\begin{align*}
M_3:=C_r2^{-\frac{2}{pq-1}-r}(pq)^{-\frac{(4q+2)pq}{(pq-1)^2}}M_2^{\frac{1}{pq-1}}>0,
\end{align*}
we can rewrite the last lower bound estimate as
\begin{align}\label{Lower-Bound-V-2}
\ml{V}(t)\geqslant\exp\left[(pq)^j\ln \left( M_3\varepsilon\,t^{\frac{2}{pq-1}+r}\right) \right](t-L)^{-\frac{2}{pq-1}}
\end{align}
for any $j\geqslant j_0$.

Let us consider $\varepsilon_0=\varepsilon_0(u_0,u_1,v_0,v_1,b,m^2,p,q)>0$ such that
\begin{align*}
\varepsilon_0\leqslant (2L)^{-\frac{2}{pq-1}-r}M_3^{-1}.
\end{align*}
Hence, for any $\varepsilon\in(0,\varepsilon_0]$ and for $t>(M_3\varepsilon)^{-\frac{pq-1}{rpq+2-r}}$, we have $t\geqslant 2L$, and the argument of the logarithmic function in \eqref{Lower-Bound-V-2} is greater than $1$. So, letting $j\to+\infty$ in \eqref{Lower-Bound-V-2}, we claim that the lower bound for the functional $\ml{V}(t)$ blows up in finite time, moreover, the lifespan of local solution can be estimated from the above as follows:
\begin{align*}
T(\varepsilon)\lesssim \varepsilon^{-\frac{pq-1}{rpq+2-r}}.
\end{align*}
We completed the proof of Theorem \ref{Thm-Blow-up} for the damped Klein-Gordon case.

\section{Proof of Theorem \ref{Thm-Blow-up massless}}
\label{Section-Blow-up-damping}\setcounter{equation}{0}
\subsection{Iteration frames and first lower bounds}
\hspace{5mm}Let us define the time-dependent functionals $\ml{U}(t)$ and $\ml{V}(t)$ as those in \eqref{Functionals}. By following the same procedure as Section \ref{Section-Blow-up-damping-mass} with $m^2=0$, we are able to derive the iteration frames
\begin{align}
\ml{U}(t)&\geqslant \mathrm{e}^{-bt}\int_0^t\mathrm{e}^{bs}\int_0^s|\ml{V}(\tau)|^p\,\mathrm{d}\tau\,\mathrm{d}s,\label{Frame-U-2}\\
\ml{V}(t)&\geqslant \int_0^t\int_0^s|\ml{U}(\tau)|^q\,\mathrm{d}\tau\,\mathrm{d}s,\label{Frame-V-2}
\end{align}
and the first lower bounds
\begin{align}
\ml{U}(t)&\geqslant \ml{U}(0)+\frac{1}{b}\big(1-\mathrm{e}^{-bt}\big)\ml{U}'(0)\geqslant C_2\varepsilon, \label{star}\\
\ml{V}(t)&\geqslant \ml{V}(0)+\ml{V}'(0)\,t\geqslant C_r\varepsilon \,t^r, \label{star 2}
\end{align}
for $t\geqslant\frac{1}{b}$, where either $r=0$ if $v_1\equiv 0$ or $r=1$ otherwise (here the non-negativity assumptions on initial data is used). Comparing \eqref{Frame-U-2} with \eqref{Frame-U}, we see that in the massless case only one unbounded exponential multiplier appears. Moreover, there is  a bounded lower bound for $\ml{U}(t)$ instead of an exponential lower bound as in  \eqref{Lower-Bound-U-1}. These phenomena will greatly affect  the lifespan of the solutions when $m^2=0$.

\subsection{Iteration argument}
\hspace{5mm} Let us consider the sequence $\{\ell_k\}_{k\in\mathbb{N}_0}$, defined as follows:
\begin{align*}
\ell_0:= \frac{1}{b} \ \ \mbox{and}\ \  \ell_k:=1+(pq)^{-k}\ \ \mbox{for any}\ \ k\geqslant 1.
\end{align*} Moreover, we employ the sequence $\{L_j\}_{j\in\mathbb{N}_0}$ which is formally defined as in \eqref{Defn-Lj}. In comparison to the previous case the difference lies in the above defined sequence $\{\ell_k\}_{k\in\mathbb{N}_0}$ since in this case we have to apply a 1-step slicing procedure.

We will determine sequences of lower bounds for $\ml{U}(t)$ and $\ml{V}(t)$, respectively, in an iterative way. More precisely, we will prove that
\begin{align}
\ml{U}(t)&\geqslant H_j(t- L_j)^{a_j}\ \ \mbox{for any}\ \ t\geqslant L_j,\label{Seq-U-2}\\
\ml{V}(t)&\geqslant K_j(t- L_j)^{b_j}\ \ \mbox{for any}\ \ t\geqslant L_j,\label{Seq-V-2}
\end{align}
where $\{H_j\}_{j\in\mb{N}_0}$, $\{K_j\}_{j\in\mb{N}_0}$, $\{a_j\}_{j\in\mb{N}_0}$ and $\{b_j\}_{j\in\mb{N}_0}$ are sequences of non-negative real numbers to be determined iteratively. As in the previous result, we will show the validity of \eqref{Seq-U-2} and \eqref{Seq-V-2} by induction on $j$. 

The base case follows by \eqref{star} and \eqref{star 2} taking $H_0:=C_2 \varepsilon$, $K_0:= C_r \varepsilon$, $a_0:=0$ and $b_0:= r$ (here $r$ has the same meaning as before). 

Let us begin with the inductive step for \eqref{Seq-U-2}. Plugging \eqref{Seq-U-2} into \eqref{Frame-V-2}, for $t\geqslant L_j$ we have
\begin{align*}
\ml{V}(t) & \geqslant \int_{L_j}^t \int_{L_j}^s (\ml{U}(\tau))^q\, \mathrm{d}\tau \, \mathrm{d}s \\ & \geqslant H_j^q \int_{L_j}^t \int_{L_j}^s (\tau-L_j)^{q a_j}\, \mathrm{d}\tau \, \mathrm{d}s =\frac{H_j^q}{(q a_j+1)(q a_j+2)} (t-L_j)^{q a_j+2}.
\end{align*}
Hence, using the last lower bound for $\ml{V}$ in \eqref{Frame-U-2}, for $t\geqslant L_{j+1}$ we have
\begin{align*}
\ml{U}(t) & \geqslant \mathrm{e}^{-bt} \int_{L_j}^t  \mathrm{e}^{bs} \int_{L_j}^s (\ml{V}(\tau))^p\, \mathrm{d}\tau \, \mathrm{d}s \\ 
& \geqslant \frac{H_j^{pq}}{(q a_j+1)^p(q a_j+2)^p} \mathrm{e}^{-bt} \int_{L_j}^t  \mathrm{e}^{bs} \int_{L_j}^s (\tau-L_j)^{pq a_j+2p}\, \mathrm{d}\tau \, \mathrm{d}s \\ 
& \geqslant \frac{H_j^{pq}}{(q a_j+1)^p(q a_j+2)^p(pq a_j+2p+1)} \mathrm{e}^{-bt} \int_{L_j}^t  \mathrm{e}^{bs}  (s-L_j)^{pq a_j+2p+1}\, \mathrm{d}s\\ 
& \geqslant \frac{H_j^{pq}}{(q a_j+1)^p(q a_j+2)^p(pq a_j+2p+1)} \mathrm{e}^{-bt} \int_{\tfrac{t}{\ell_{j+1}}}^t  \mathrm{e}^{bs}  (s-L_j)^{pq a_j+2p+1}\, \mathrm{d}s\\ 
& \geqslant \frac{H_j^{pq} \ell_{j+1}^{-(pq a_j+2p+1)}}{b(q a_j+1)^p(q a_j+2)^p(pq a_j+2p+1)}  (t-L_{j+1})^{pq a_j+2p+1} \left(1-\mathrm{e}^{-b(\ell_{j+1}-1)\frac{t}{\ell_{j+1}}}\right) \\
& \geqslant \frac{H_j^{pq} \ell_{j+1}^{-(pq a_j+2p+1)}}{b(pq a_j+2p+1)^{2p+1}}  (t-L_{j+1})^{pq a_j+2p+1} \left(1-\mathrm{e}^{-(\ell_{j+1}-1)}\right),
\end{align*} where we shrank the domain of integration in the $s$-integral from $[L_j,t]$ to $[\frac{t}{\ell_{j+1}},t]$ and we used $b\frac{t}{\ell_{j+1}}\geqslant b L_j\geqslant 1$. Since
\begin{align}
1-\mathrm{e}^{-(\ell_{j+1}-1)} &\geqslant (\ell_{j+1}-1)\big(1-\tfrac{1}{2}(\ell_{j+1}-1)\big)  \notag \\
&\geqslant (pq)^{-2(j+1)}\big((pq)^{j+1}-\tfrac{1}{2}\big) \notag  \\
&\geqslant(pq)^{-2(j+1)}\big(pq-\tfrac{1}{2}\big), \label{estimate exponential term}
\end{align} we proved \eqref{Seq-U-2} for $j+1$ provided that
\begin{align}
a_{j+1}  := pq a_j+2p+1 \ \ \mbox{and}\ \ H_{j+1} & := \frac{ (pq-\tfrac{1}{2})  \ell_{j+1}^{-a_{j+1}}}{b \,a_{j+1}^{2p+1} (pq)^{2(j+1)}} H_j^{pq}.  \label{(II)}
\end{align} We prove now the inductive step for \eqref{Seq-V-2}. Combining \eqref{Seq-V-2} for $j\geqslant 0$ and \eqref{Frame-U-2}, for $t\geqslant L_{j+1}$ we find
\begin{align*}
\ml{U}(t) & \geqslant  \mathrm{e}^{-bt}\int_{L_{j}}^t\mathrm{e}^{bs}\int_{L_{j}}^s (\ml{V}(\tau))^p\,\mathrm{d}\tau\,\mathrm{d}s \\ & \geqslant K_j^p  \,\mathrm{e}^{-bt}\int_{L_{j}}^t\mathrm{e}^{bs}\int_{L_{j}}^s(\tau-L_j)^{p b_j}\,\mathrm{d}\tau\,\mathrm{d}s \\
& \geqslant \frac{K_j^p}{p b_j+1} \, \mathrm{e}^{-bt}\int_{L_{j}}^t\mathrm{e}^{bs} (s-L_j)^{p b_j+1}\,\mathrm{d}s \\
& \geqslant \frac{K_j^p}{p b_j+1} \, \mathrm{e}^{-bt}\int_{\tfrac{t}{\ell_{j+1}}}^t\mathrm{e}^{bs} (s-L_j)^{p b_j+1}\,\mathrm{d}s \\
& \geqslant \frac{K_j^p \ell_{j+1}^{-(p b_j+1)}}{b(p b_j+1)} \, (t-L_{j+1})^{p b_j+1} \left(1-\mathrm{e}^{-b(\ell_{j+1}-1)\frac{t}{\ell_{j+1}}}\right)\\
& \geqslant \frac{ (pq-\tfrac{1}{2})  \ell_{j+1}^{-(p b_j+1)} K_j^p }{b(p b_j+1)  (pq)^{2(j+1)}} \, (t-L_{j+1})^{p b_j+1},
\end{align*} where we shrank the domain as before and we used the estimate  \eqref{estimate exponential term}.
Applying the previous estimate in \eqref{Frame-V-2}, for $t\geqslant L_{j+1}$ it follows that
\begin{align*}
\ml{V}(t) & \geqslant \int_{L_{j+1}}^t\int_{L_{j+1}}^s(\ml{U}(\tau))^q\,\mathrm{d}\tau\,\mathrm{d}s \\
& \geqslant \frac{ (pq-\tfrac{1}{2})^q  \ell_{j+1}^{-q(p b_j+1)} K_j^{pq} }{b^q(p b_j+1)^q  (pq)^{2q(j+1)}}  \int_{L_{j+1}}^t\int_{L_{j+1}}^s  (\tau-L_{j+1})^{pq b_j+q} \,\mathrm{d}\tau\,\mathrm{d}s \\
& \geqslant \frac{ (pq-\tfrac{1}{2})^q  \ell_{j+1}^{-q(p b_j+1)} K_j^{pq} }{b^q(p b_j+1)^q (pq b_j+q+1) (pq b_j+q+2) (pq)^{2q(j+1)}}   (t-L_{j+1})^{pq b_j+q+2} \\
& \geqslant \frac{ (pq-\tfrac{1}{2})^q  \ell_{j+1}^{-(pq b_j+q+2)} K_j^{pq} }{b^q (pq b_j+q+2)^{q+2} (pq)^{2q(j+1)}}   (t-L_{j+1})^{pq b_j+q+2}.
\end{align*} Thus, we proved \eqref{Seq-V-2} for $j+1$ provided that
\begin{align}
b_{j+1}  := pq b_j+q+2\ \ \mbox{and}\ \ K_{j+1} & := \frac{ (pq-\tfrac{1}{2})^q  \ell_{j+1}^{- b_{j+1} }  }{b^q b_{j+1}^{q+2} (pq)^{2q(j+1)}}K_j^{pq}. \label{(IV)}
\end{align}

\subsection{Upper bound estimates for the lifespan}
\hspace{5mm}By using the recursive relation $a_j=2p+1+pq a_{j-1}$ in an iterative way, we obtain
\begin{align}
a_j &= (2p+1)(1+pq)+(pq)^2 a_{j-2} = (2p+1)\big(1+pq+(pq)^2\big)+(pq)^3 a_{j-3} \notag \\& = \cdots= (2p+1) \sum_{k=0}^{j-1} (pq)^k+ (pq)^j a_0 = \left(a_0+\frac{2p+1}{pq-1}\right)(pq)^j-\frac{2p+1}{pq-1}. \label{(V)}
\end{align} Similarly,
\begin{align}
b_j= \left(b_0+\frac{q+2}{pq-1}\right)(pq)^j-\frac{q+2}{pq-1}. \label{(VI)}
\end{align}
Hence, from \eqref{(V)} and \eqref{(VI)} we get
\begin{align}
a_j^{2p+1} & \leqslant \left(a_0+\frac{2p+1}{pq-1}\right)^{2p+1} (pq)^{(2p+1)j}, \label{(VII)} \\
b_j^{q+2} & \leqslant \left(b_0+\frac{q+2}{pq-1}\right)^{q+2} (pq)^{(q+2)j}. \label{(VIII)}
\end{align}
Furthermore, we remark that
\begin{align*}
\lim_{j\to +\infty} \ell_j^{a_j} & = \lim_{j\to +\infty} \exp\left(a_j \ln \ell_j\right) \\
& = \lim_{j\to +\infty} \exp\left( \left(a_0+\tfrac{2p+1}{pq-1}\right)(pq)^j \ln \left(1+(pq)^{-j}\right)\right) = \exp\left(a_0+\tfrac{2p+1}{pq-1}\right),
\end{align*} and, analogously,
\begin{align*}
\lim_{j\to +\infty} \ell_j^{b_j} &  = \exp\left(b_0+\tfrac{q+2}{pq-1}\right).
\end{align*}
Therefore, there exists $M_4=M_4(p,q,v_0,v_1)>0$ such that $\ell_{j}^{-a_j},\ell_{j}^{-b_j}\leqslant M_4$ for any $j\in\mathbb{N}_0$.

In conclusion, combining \eqref{(II)}, \eqref{(IV)}, \eqref{(VII)} and \eqref{(VIII)} with the (uniform in $j$) lower bound for  $\ell_{j}^{-a_j},\ell_{j}^{-b_j}$, we obtain
\begin{align*}
H_j & \geqslant H (pq)^{-(2p+3)j} H_{j-1}^{pq}, \\
K_j & \geqslant K (pq)^{-(3q+2)j} K_{j-1}^{pq}, 
\end{align*} where $H:= \frac{M_4}{b}(pq-\frac{1}{2})\left(a_0+\frac{2p+1}{pq-1}\right)^{-(2p+1)}$ and $K:= \frac{M_4}{b^q}(pq-\frac{1}{2})^q\left(b_0+\frac{q+2}{pq-1}\right)^{-(q+2)}$.
We repeat similar computations to those done in Subsection \ref{Subsection 3.3}: applying the logarithmic function to both sides of the previous inequalities for $H_j$ and $K_j$, we find
\begin{align*}
\ln H_j & \geqslant (pq) \ln H_{j-1}-(2p+3)j\ln(pq)+\ln H, \\
\ln K_j & \geqslant (pq) \ln K_{j-1}-(3q+2)j\ln(pq)+\ln K. 
\end{align*}
Let $j_1,j_2\in \mathbb{N}_0$ be the smallest integers satisfying $j_1\geqslant \frac{\ln H}{(2p+3)\ln(pq)}-\frac{pq}{pq-1}$ and $j_2\geqslant \frac{\ln K}{(3q+2)\ln(pq)}-\frac{pq}{pq-1}$, respectively. Then,
\begin{align}
\ln H_j & \geqslant (pq)^j \left(\ln H_0-\frac{(2p+3)pq}{(pq-1)^2}\ln(pq)+\frac{\ln H}{pq-1}\right)\ \  \mbox{for any } j\geqslant j_1, \label{(IX)} \\
\ln K_j & \geqslant (pq)^j \left(\ln K_0-\frac{(3q+2)pq}{(pq-1)^2}\ln(pq)+\frac{\ln K}{pq-1}\right)\ \  \mbox{for any } j\geqslant j_2. \label{(X)} 
\end{align}
Since $L_j\uparrow L:= \prod_{k=0}^{+\infty}\ell_k$ as $j\to+\infty$, from \eqref{Seq-U-2}, \eqref{(V)} and \eqref{(IX)} for $j\geqslant j_1$ and $t\geqslant L$ we have
\begin{align*}
\ml{U}(t) & \geqslant H_j (t-L)^{a_j} \\
& \geqslant \exp \left[ (pq)^j \left(\ln H_0 -\frac{(2p+3)pq}{(pq-1)^2}\ln(pq)+\frac{\ln H}{pq-1}+\left(a_0+\frac{2p+1}{pq-1}\right)\ln(t-L)\right)\right] (t-L)^{-\frac{2p+1}{pq-1}}.
\end{align*} Then, for $t\geqslant 2L$, so that $\ln(t-L)\geqslant \ln t-\ln 2$, and for $j\geqslant j_1$, it results
\begin{align} \label{(XI)}
\ml{U}(t) & \geqslant \exp \left[ (pq)^j \ln\left(M_5 \varepsilon\, t^{a_0+\frac{2p+1}{pq-1}}\right)\right] (t-L)^{-\frac{2p+1}{pq-1}},
\end{align} where $M_5:= C_2 (pq)^{-\frac{(2p+3)pq}{(pq-1)^2}}2^{-(a_0+\frac{2p+1}{pq-1})}H^{\frac{1}{pq-1}}$. 

Similarly, by \eqref{Seq-V-2}, \eqref{(VI)} and \eqref{(X)}, for $t\geqslant 2L$ and for $j\geqslant j_2$ it holds
\begin{align} \label{(XII)}
\ml{V}(t) & \geqslant \exp \left[ (pq)^j \ln\left(M_6 \varepsilon\, t^{b_0+\frac{q+2}{pq-1}}\right)\right] (t-L)^{-\frac{q+2}{pq-1}},
\end{align} where $M_6:= C_r (pq)^{-\frac{(3q+2)pq}{(pq-1)^2}}2^{-(b_0+\frac{q+2}{pq-1})}K^{\frac{1}{pq-1}}$.

Finally, we fix $\varepsilon_0=\varepsilon_0(u_0,u_1,v_0,v_1,b,p,q)>0$ such that $$\varepsilon_0\leqslant \min\left\{M_5^{-1}(2L)^{-\left(a_0+\frac{2p+1}{pq-1}\right)},M_6^{-1}(2L)^{-\left(b_0+\frac{q+2}{pq-1}\right)}\right\}.$$
Then, for any $\varepsilon\in (0,\varepsilon_0]$, whereas for $t>(M_5\varepsilon)^{-\frac{pq-1}{2p+1}}$ the right-hand side of \eqref{(XI)} blows up as $j\to+\infty$, for $t>(M_6\varepsilon)^{-\frac{pq-1}{rpq+2-r+q}}$ the right-hand side of \eqref{(XII)} blows up as $j\to+\infty$. So, we proved that $(\ml{U},\ml{V})$ blows up in finite time and for any $\varepsilon\in(0,\varepsilon_0]$ we showed the upper bound estimates for the lifespan
\begin{align*}
T(\varepsilon) \lesssim \min \left\{\varepsilon^{-\frac{pq-1}{2p+1}}, \varepsilon^{-\frac{pq-1}{rpq+2-r+q}} \right\} = \varepsilon^{-\frac{pq-1}{\max\{2p+1,rpq +2-r+q\}}},
\end{align*} where either $r=0$ if $v_1\equiv 0$ or $r=1$ otherwise.

\section{Local in-time solution  for the wave system with lower order terms}\setcounter{equation}{0}
\subsection{Preliminary results}
\hspace{5mm}Let us recall the notion of mild solutions for the weakly coupled system \eqref{Eq-Main-Systems}. Let
\begin{align*}
N_1v:=u^{\lin}+J_1v\ \ \mbox{and} \ \ 
N_2u:=v^{\lin}+J_2u,
\end{align*}
where
\begin{align}
	u^{\lin}(t,x)&:=\varepsilon u_0(x)\ast_{(x)}E_0(t,x;b,m^2)+\varepsilon u_1(x)\ast_{(x)}E_1(t,x;b,m^2),\label{ulin}\\
	v^{\lin}(t,x)&:=\varepsilon v_0(x)\ast_{(x)}E_0(t,x;0,0)+\varepsilon v_1(x)\ast_{(x)}E_1(t,x;0,0),\label{vlin}
\end{align}
and
\begin{align*}
	J_1 v(t,x)&:=\int_0^t|v(s,x)|^p\ast_{(x)}E_1(t-s,x;b,m^2)\,\mathrm{d}s, \\
	J_2 u(t,x)&:=\int_0^t|u(s,x)|^q\ast_{(x)}E_1(t-s,x;0,0)\,\mathrm{d}s.
\end{align*}
Here, for $b\geqslant 0$ and $m^2\geqslant 0$, $E_0(t,x;b,m^2)$ and $E_1(t,x;b,m^2)$ denote the fundamental solutions to the linear Cauchy problem for $\partial_t^2w-\ml{L}w+b\,\partial_tw+m^2w=0$ on a compact Lie group $\G$
with the initial data $(\delta_0,0)$ and $(0,\delta_0)$ (see \cite{Palmieri-JDE=2021,Palmieri-Pseudo=2021,Palmieri-JFAA=2022} for further details). In other words, $(u^{\lin},v^{\lin})$ is the solution to the homogeneous wave system associated with \eqref{Eq-Main-Systems}.
Motivated by Duhamel's principle, we define the operator
\begin{align*}
\ml{N}:\ (u,v)\in X(T)\to\ml{N}[u,v]:=(N_1 v,N_2 u)
\end{align*}
on a suitable evolution space $X(T):=Y_1(T)\times Y_2(T)$ whose definition will depend on whether we are in the damped Klein-Gordon case or in the massless case.

\begin{defn}\label{Defn-Mild}
	Let $b>0$ and $m^2\geqslant 0$. We say that $(u,v)$ is a mild solution to the weakly coupled system \eqref{Eq-Main-Systems} on $[0,T]$ if $(u,v)$ is a fixed point for $\ml{N}$ on $X(T)$, that is, if $(u,v)$ is a solution to the nonlinear integral system
	\begin{align*}
		\begin{cases}
			u=u^{\lin}+J_1v,\\
			v=v^{\lin}+J_2u.
		\end{cases}
	\end{align*}
\end{defn}
We point out that, in order to get the previous representation formulas, we applied the invariance by time translations for the differential operators $\partial_t^2-\ml{L}+b\,\partial_t+m^2$ and $\partial_t^2-\ml{L}$. Moreover, the identity $$\mathit{L}\big(w\ast_{(x)}E_1(t,\cdot\,;b,m^2)\big)=w\ast_{(x)}\mathit{L}\big(E_1(t,\cdot\,;b,m^2)\big)$$ for any left-invariant differential operator $\mathit{L}$ on $\G$ is used to derive the integral formulation of the problem.

We set
\begin{align*}
\Lambda:=\big\{(j,k)\in \mathbb{N}_0:\ \ 0\leqslant j+k\leqslant 1\big\}.
\end{align*}
 Next, we recall the $L^2(\G)-L^2(\G)$ estimates for the solutions to the linear Cauchy problems for the single equations. These estimates were obtained by using the group Fourier transform with respect to the spatial variable $x\in\G$  and the Plancherel identity in the framework of compact Lie groups.

\begin{prop}[\cite{Palmieri-JDE=2021,Palmieri-JFAA=2022}]\label{Prop-Palmieri-1} Let us assume $(u_0,u_1)\in H^1_{\ml{L}}(\G)\times L^2(\G)$ and $b>0$ satisfying $b^2\geqslant 4m^2$. Then, the mild solution
	\begin{align*}
	u^{\lin}\in \ml{C}\big([0,+\infty), H^1_{\ml{L}}(\G)\big)\cap \ml{C}^1\big([0,+\infty), L^2(\G)\big)
	\end{align*} defined in \eqref{ulin} satisfies	the following $L^2(\G)-L^2(\G)$ estimates:
	\begin{align*}
		\|\partial_t^j(-\ml{L})^{k/2}u^{\lin}(t,\cdot)\|_{L^2(\G)}&\lesssim \begin{cases}
			\varepsilon	d_{b,m^2}(t)\|(u_0,u_1)\|_{H^{j+k}_{\ml{L}}(\G)\times L^2(\G)}&\mbox{if}\ \ m^2>0,\\
			\varepsilon (1+t)^{-\frac{2j+k}{2}}\|(u_0,u_1)\|_{H^{j+k}_{\ml{L}}(\G)\times L^2(\G)}&\mbox{if}\ \ m^2=0,
		\end{cases}
	\end{align*}
	for any $(j,k)\in\Lambda$,
	 where
	\begin{align*}
		d_{b,m^2}(t):=\begin{cases}
			\mathrm{e}^{-\frac{1}{2}(b-\sqrt{b^2-4m^2}\,)t}&\mbox{if}\ \ b^2>4m^2,\\
			(1+t)\,\mathrm{e}^{-\frac{b}{2}t}&\mbox{if}\ \ b^2=4m^2.
		\end{cases}
	\end{align*}
\end{prop}

\begin{prop}[\cite{Palmieri-Pseudo=2021}]\label{Prop-Palmieri-2} Let us assume $(v_0,v_1)\in H^1_{\ml{L}}(\G)\times L^2(\G)$. Then, the mild solution 
\begin{align*}
v^{\lin}\in \ml{C}\big([0,+\infty), H^1_{\ml{L}}(\G)\big)\cap \ml{C}^1\big([0,+\infty), L^2(\G)\big)
\end{align*}
 defined in  \eqref{vlin} satisfies	the following $L^2(\G)-L^2(\G)$ estimates:
	\begin{align*}
		\|\partial_t^j(-\ml{L})^{k/2}v^{\lin}(t,\cdot)\|_{L^2(\G)}&\lesssim \varepsilon[a(t)]^{1-(j+k)}\|(v_0,v_1)\|_{H^{j+k}_{\ml{L}}(\G)\times L^2(\G)},
	\end{align*}
for any $(j,k)\in\Lambda$, where
	\begin{align*}
		a(t):=\begin{cases}
			1+t&\mbox{if} \  \ v_1\neq0,\\
			1&\mbox{if}\ \ v_1=0.
		\end{cases}
	\end{align*}
\end{prop}

A useful and fundamental tool  to estimate power nonlinearities is the next Gagliardo-Nirenberg type inequality, whose proof can be found in \cite{Ruzhansky-Yess=2024} in a more general setting (see also \cite[Corollary 2.3]{Palmieri-JDE=2021}).
\begin{lemma}\label{Lemma-GN-inequality}
	Let $\G$ be a connected unimodular Lie group with topological dimension $n\geqslant 3$. For any $\gamma\geqslant 2$ such that $\gamma\leqslant\frac{2n}{n-2}$, the following Gagliardo-Nirenberg type inequality:
	\begin{align*}
		\|f\|_{L^{\gamma}(\G)}\lesssim\|f\|_{H^1_{\ml{L}}(\G)}^{\theta(n,\gamma)}\|f\|_{L^2(\G)}^{1-\theta(n,\gamma)}
	\end{align*}
	holds for any $f\in H^1_{\ml{L}}(\G)$, where $\theta(n,\gamma):=n(\frac{1}{2}-\frac{1}{\gamma})$.
\end{lemma}

\subsection{Proof of Theorem \ref{Thm-Local}}
\hspace{5mm}Motivated by the $L^2(\G)-L^2(\G)$ estimates in Propositions \ref{Prop-Palmieri-1} and \ref{Prop-Palmieri-2}, we introduce
the evolution spaces 
\begin{align}\label{Yspace}
Y_{j}(T):=\ml{C}\big([0,T], H^1_{\ml{L}}(\G)\big)\cap \ml{C}^1\big([0,T], L^2(\G)\big) \ \ \mbox{for}\ \ j\in\{1,2\}
\end{align}
 endowed with the norms
\begin{align}
	\|u\|_{Y_1(T)}&:=\sup\limits_{t\in[0,T]}\left((1+t)^{\lambda_0}\sum\limits_{(j,k)\in\Lambda}\|\partial_t^j(-\ml{L})^{k/2}u(t,\cdot)\|_{L^2(\G)}\right),\notag\\
	\|v\|_{Y_2(T)}&:=\sup\limits_{t\in[0,T]}\left(\,\sum\limits_{(j,k)\in\Lambda}[a(t)]^{j+k-1}\|\partial_t^j(-\ml{L})^{k/2}v(t,\cdot)\|_{L^2(\G)}\right),\label{Norm-Y1-Y2}
\end{align}
where the parameter $\lambda_0$ is defined by
\begin{align*}
	\lambda_0:=\begin{cases}
		\displaystyle{\frac{2p-pq-1}{pq-1}}&\mbox{if}\ \ v_1\neq0,\\[1em]
		\displaystyle{\frac{2p-2}{pq-1}}&\mbox{if}\ \ v_1=0.
	\end{cases}
\end{align*} 
Thanks to our definitions of the evolution spaces, we may estimate the power nonlinearities in $L^2(\G)$ by using the Gagliardo-Nirenberg type inequality. Consequently,
\begin{align}
	\|\,|v(s,\cdot)|^p\|_{L^2(\G)}&\lesssim \|v(s,\cdot)\|_{H^1_{\ml{L}}(\G)}^{p\theta(n,2p)}\|v(s,\cdot)\|_{L^2(\G)}^{p[1-\theta(n,2p)]}\lesssim [a(s)]^p\|v\|_{Y_2(s)}^p,\label{nonlinear-1}\\
	\|\,|u(s,\cdot)|^q\|_{L^2(\G)}&\lesssim \|u(s,\cdot)\|_{H^1_{\ml{L}}(\G)}^{q\theta(n,2q)}\|u(s,\cdot)\|_{L^2(\G)}^{q[1-\theta(n,2q)]}\lesssim(1+s)^{-\lambda_0 q}\|u\|_{Y_1(s)}^q. \label{nonlinear-2}
\end{align}
 Let us stress that the  assumptions $p,q\leqslant\frac{n}{n-2}$ and $n\geqslant 3$ are necessary for the previous estimates.

Let us first derive a local in-time existence result for mild solutions to \eqref{Eq-Main-Systems}.  By Propositions \ref{Prop-Palmieri-1} and \ref{Prop-Palmieri-2}, it follows immediately that
\begin{align}
	\|u^{\lin}\|_{Y_1(T)}&\lesssim\varepsilon\|(u_0,u_1)\|_{H^1_{\ml{L}}(\G)\times L^2(\G)},\label{L1}\\
	\|v^{\lin}\|_{Y_2(T)}&\lesssim\varepsilon\|(v_0,v_1)\|_{H^1_{\ml{L}}(\G)\times L^2(\G)},\label{L2}
\end{align}
because of $(1+t)^{\lambda_0}d_{b,m^2}(t)\lesssim 1$. On the other hand, thanks to the invariance by time translations of the corresponding linear Cauchy problems, by using \eqref{nonlinear-1}, we get 
\begin{align*}
	&(1+t)^{\lambda_0}\|\partial_t^j(-\ml{L})^{k/2}J_1 v(t,\cdot)\|_{L^2(\G)}\\
	&\qquad\lesssim(1+t)^{\lambda_0}\int_0^td_{b,m^2}(t-s)\|\,|v(s,\cdot)|^p\|_{L^2(\G)}\,\mathrm{d}s\\
	&\qquad\lesssim \left((1+t)^{\lambda_0}d_{b,m^2}(t/2)\int_0^{t/2}[a(s)]^p\,\mathrm{d}s+(1+t)^{\lambda_0}[a(t)]^p\int_{t/2}^td_{b,m^2}(t-s)\,\mathrm{d}s\right)\|v\|_{Y_2(t)}^p\\
	&\qquad\lesssim\left((1+t)^{\lambda_0+1}d_{b,m^2}(t/2)[a(t)]^{p}+(1+t)^{\lambda_0}[a(t)]^p\right)\|v\|_{Y_2(t)}^p\\
	&\qquad\lesssim (1+t)^{\lambda_0}[a(t)]^p\|v\|_{Y_2(t)}^p
\end{align*}
for any $(j,k)\in\Lambda$, where we applied the estimates
\begin{align*}
d_{b,m^2}(t-s)\leqslant
\begin{cases}
d_{b,m^2}(t/2)&\mbox{if}\ \ b^2>4m^2,\\
2d_{b,m^2}(t/2)&\mbox{if}\ \ b^2=4m^2,
\end{cases}
\end{align*}
for $s\in[0,t/2]$, and the integrability of $d_{b,m^2}(t-s)$ for $s\in[t/2,t]$ due to
\begin{align*}
\int_{t/2}^td_{b,m^2}(t-s)\,\mathrm{d}s=\int_0^{t/2}d_{b,m^2}(s)\,\mathrm{d}s\begin{cases}
\lesssim 1-d_{b,m^2}(t/2)\lesssim 1&\mbox{if}\ \ b^2>4m^2,\\
=\frac{2}{b}\left[(1+\frac{2}{b})-(1+\frac{2}{b}+\frac{t}{2})\,\mathrm{e}^{-\frac{b}{4}t}\right]\lesssim 1&\mbox{if}\ \ b^2=4m^2.
\end{cases}
\end{align*}
Thanks to the choice of $\lambda_0$, it holds
\begin{align*}
	1-\lambda_0q=\begin{cases}
		\displaystyle{\frac{(pq+1)(q-1)}{pq-1}> 0}&\mbox{if}\ \ v_1\neq0,\\[1em]
		\displaystyle{\frac{-pq-1+2q}{pq-1}>0}&\mbox{if}\ \ v_1=0,
	\end{cases}
\end{align*}
where the further condition $p<2-\frac{1}{q}$ is assumed when $v_1=0$.
 Analogously, from \eqref{nonlinear-2} we derive
\begin{align*}
	&[a(t)]^{j+k-1}\|\partial_t^j(-\ml{L})^{k/2}J_2 u(t,\cdot)\|_{L^2(\G)}\\
	&\qquad\lesssim[a(t)]^{j+k-1}\int_0^t(1+t-s)^{1-(j+k)}\|\,|u(s,\cdot)|^q\|_{L^2(\G)}\,\mathrm{d}s\\
	&\qquad\lesssim\left([a(t)]^{j+k-1}(1+t)^{1-(j+k)}\int_0^{t/2}(1+s)^{-\lambda_0 q}\,\mathrm{d}s\right.\\
	&\qquad\qquad\left.+[a(t)]^{j+k-1}(1+t)^{-\lambda_0 q}\int_{t/2}^t(1+t-s)^{1-(j+k)}\,\mathrm{d}s\right)\|u\|_{Y_1(t)}^q\\
	&\qquad\lesssim (1+t)^{2-(j+k)-\lambda_0q}[a(t)]^{j+k-1}\|u\|_{Y_1(t)}^q
\end{align*}
for any $(j,k)\in\Lambda$. Summarizing the last estimates, we have proved that
\begin{align}
	\|N_1 v\|_{Y_1(T)}&\lesssim \varepsilon\|(u_0,u_1)\|_{H^1_{\ml{L}}(\G)\times L^2(\G)} +
		(1+T)^{\lambda_0}[a(T)]^p\|v\|_{Y_2(T)}^p,\label{N1u-v}\\
	\|N_2 u\|_{Y_2(T)}&\lesssim \varepsilon\|(v_0,v_1)\|_{H^1_{\ml{L}}(\G)\times L^2(\G)} +(1+T)^{2-\lambda_0q}[a(T)]^{-1}\|u\|_{Y_1(T)}^q.\label{N2v-u}
\end{align}
By using the mean value theorem for the $p$ and $q$-powers and H\"older's inequality, for any $u,\tilde{u}\in Y_1(T)$ and $v,\tilde{v}\in Y_2(T)$, one can derive the following inequalities:
\begin{align*}
	\|N_1 v-N_1\tilde{v}\|_{Y_1(T)}&\lesssim (1+T)^{\lambda_0}[a(T)]^p\|v-\tilde{v}\|_{Y_2(T)}\left(\|v\|_{Y_2(T)}^{p-1}+\|\tilde{v}\|_{Y_2(T)}^{p-1}\right),\\
	\|N_2 u-N_2\tilde{u}\|_{Y_2(T)}&\lesssim (1+T)^{2-\lambda_0q}[a(T)]^{-1}\|u-\tilde{u}\|_{Y_1(T)}\left(\|u\|_{Y_1(T)}^{q-1}+\|\tilde{u}\|_{Y_1(T)}^{q-1}\right).
\end{align*}
Consequently, defining $\| (u,v)\|_{X(T)}:= \max\{\|u\|_{Y_1(T)}, \|v\|_{Y_2(T)}\}$ for any $(u,v)\in X(T)$, then we can summarize the previous inequalities as follows: 
\begin{align}
	& \|\ml{N}[u,v]\|_{X(T)} \lesssim\varepsilon\left(\|(u_0,u_1)\|_{H^1_{\ml{L}}(\G)\times L^2(\G)}+ \|(v_0,v_1)\|_{H^1_{\ml{L}}(\G)\times L^2(\G)}\right)+\sum\limits_{\gamma\in\{p,q\}}\|(u,v)\|_{X(T)}^{\gamma}, \label{contraction 1}\\
	& \|\ml{N}[u,v]-\ml{N}[\tilde{u},\tilde{v}]\|_{X(T)} \lesssim \|(u,v)-(\tilde{u},\tilde{v})\|_{X(T)}\left(\, \sum\limits_{\gamma\in\{p,q\}}\|(u,v)\|_{X(T)}^{\gamma-1}+\|(\tilde{u},\tilde{v})\|_{X(T)}^{\gamma-1}\right), \label{contraction 2}
\end{align}
for 
$T>0$ in a bounded set, because $(1+T)^{\lambda_0}[a(T)]^p$ and $(1+T)^{2-\lambda_0q}[a(T)]^{-1}$ are regular with respect to $T$. 

By using a standard approach, provided that $R \approx\varepsilon$ 	(with independent of $T$ multiplicative constants), it is known that, on the one hand, from \eqref{contraction 1} we have that $\ml{N}$ maps 
\begin{align*}
\ml{B}_R:=\big\{(u,v)\in X(T):\  \|(u,v)\|_{X(T)}\leqslant R\big\}
\end{align*}
into itself, and, on the other hand, from \eqref{contraction 2} we have that $\ml{N}$ is a contraction. Let us stress that, although here we are working with a small-data solution ($R \approx\varepsilon$), with minor modifications it is possible to have the local existence for small time, regardless of the size of the data (i.e. without any further condition on $R$). Here, we can additionally require the smallness condition $R \approx\varepsilon$ (without loss of generality) since we are interested in establishing the lower bound estimates for the lifespan (which are meaningful just for $\varepsilon \to 0^+$).

Let $(u,v)\in  X(T)$ be the unique fixed point of $\ml{N}$ in $\ml{B}_R$ with $R\approx \varepsilon$. Combining \eqref{N1u-v} and \eqref{N2v-u} we have
	\begin{align*}
	\|u\|_{Y_1(T)}&\lesssim \varepsilon+(1+T)^{\lambda_0}[a(T)]^p\varepsilon^p+(1+T)^{2p-\lambda_0(pq-1)}\|u\|_{Y_1(T)}^{pq},\\
	\|v\|_{Y_2(T)}&\lesssim \varepsilon+(1+T)^{2-\lambda_0 q}[a(T)]^{-1}\varepsilon^q+ (1+T)^{2}[a(T)]^{pq-1}\|v\|_{Y_2(T)}^{pq},
\end{align*}
where the unexpressed multiplicative constants are independent of $T$ and $\varepsilon$. 
 In order to guarantee that $\ml{N}$ maps $\ml{B}_R$ (for $R\approx \varepsilon$) into itself, the conditions $\|u\|_{Y_1(T)}\lesssim \varepsilon$ and $\|v\|_{Y_2(T)}\lesssim \varepsilon $ are necessary. Therefore, we have to impose the following constraints on $T$:
 	\begin{align*}
 \begin{cases}
  	(1+T)^{\lambda_0}[a(T)]^p\varepsilon^p\lesssim \varepsilon,\\
 (1+T)^{2p-\lambda_0(pq-1)}\varepsilon^{pq}\lesssim\varepsilon,\\
 (1+T)^{2-\lambda_0q}[a(T)]^{-1}\varepsilon^q\lesssim \varepsilon,\\ 
 (1+T)^{2}[a(T)]^{pq-1}\varepsilon^{pq}\lesssim\varepsilon.
 \end{cases}
 \end{align*}
 Due to our definition of $\lambda_0$, these four inequalities imply the same condition on $T$, namely,
 	\begin{align*}
 	T\leqslant 1+T\lesssim \begin{cases}
 		\varepsilon^{-\frac{pq-1}{pq+1}}&\mbox{if}\ \ v_1\neq0,\\
 		\varepsilon^{-\frac{pq-1}{2}}&\mbox{if}\ \ v_1=0,
 	\end{cases}
 	\end{align*}
 	where again the unexpressed multiplicative constants are independent of $T$ and $\varepsilon$. Then, we completed the proof of the lower bound estimates for the lifespan $T(\varepsilon)$ in Theorem \ref{Thm-Local}.

\subsection{Proof of Theorem \ref{Thm-Local massless case}}
\hspace{5mm}Let us define the evolution space 
\begin{align*}
Y_1(T):=\ml{C}\big([0,T], H^1_{\ml{L}}(\G)\big)
\end{align*} and $Y_2(T)$ as in \eqref{Yspace}.
We endow $Y_2(T)$ with the same norm as in \eqref{Norm-Y1-Y2}, and $Y_1(T)$ with the following norm:
\begin{align}\label{Norm-Y1-new}
\|u\|_{Y_1(T)}:=\sup\limits_{t\in[0,T]}\left(\,\sum\limits_{k\in\{0,1\}}(1+t)^{\frac{k}{2}}\|(-\ml{L})^{k/2}u(t,\cdot)\|_{L^2(\G)}\right).
\end{align}
In other words, in the massless case we have no-loss of decay. Applying as before the Gagliardo-Nirenberg type inequality, it results
\begin{align}\label{GN-3}
\|\,|u(s,\cdot)|^q\|_{L^2(\G)}\lesssim \|u\|_{Y_1(s)}^q.
\end{align}
Clearly, the estimates \eqref{L1} and \eqref{L2} still hold due to our choices for the weights in the norms of $Y_1(T)$ and $Y_2(T)$. From Propositions \ref{Prop-Palmieri-1} and \ref{Prop-Palmieri-2}, one arrives at
	\begin{align*}
	&(1+t)^{\frac{k}{2}}\|(-\ml{L})^{k/2}J_1v(t,\cdot)\|_{L^2(\G)}\\
	&\qquad\lesssim \left(\int_0^{t/2}[a(s)]^p\,\mathrm{d}s+(1+t)^{\frac{k}{2}}[a(t)]^p\int_{t/2}^t(1+t-s)^{-\frac{k}{2}}\,\mathrm{d}s\right)\|v\|_{Y_2(t)}^p\\
	&\qquad\lesssim (1+t)[a(t)]^p\|v\|_{Y_2(t)}^p
\end{align*}
for any $k\in\{0,1\}$.  Analogously, one uses \eqref{GN-3} to derive
\begin{align*}
	&[a(t)]^{j+k-1}\|\partial_t^j(-\ml{L})^{k/2}J_2u(t,\cdot)\|_{L^2(\G)}\\
	&\qquad\lesssim\left([a(t)]^{j+k-1}(1+t)^{1-(j+k)}\int_0^{t/2}\,\mathrm{d}s+[a(t)]^{j+k-1}\int_{t/2}^t(1+t-s)^{1-(j+k)}\,\mathrm{d}s\right)\|u\|_{Y_1(t)}^q\\
	&\qquad\lesssim (1+t)^{2-(j+k)}[a(t)]^{j+k-1}\|u\|_{Y_1(t)}^q
\end{align*}
for any $(j,k)\in\Lambda$. Summarizing,
\begin{align}\label{N1u-v-2}
\|N_1 v\|_{Y_1(T)}&\lesssim \varepsilon\|(u_0,u_1)\|_{H^1_{\ml{L}}(\G)\times L^2(\G)} +
(1+T)[a(T)]^p\|v\|_{Y_2(T)}^p,\\
\|N_2 u\|_{Y_2(T)}&\lesssim\varepsilon\|(v_0,v_1)\|_{H^1_{\ml{L}}(\G)\times L^2(\G)} +(1+T)^2[a(T)]^{-1}\|u\|_{Y_1(T)}^q.\label{N1u-v-2-2}
\end{align}
The local in-time existence of mild solutions can be proved as in Theorem \ref{Thm-Local}. 

Then, the combination of \eqref{N1u-v-2} and \eqref{N1u-v-2-2} provides
	\begin{align*}
	\|u\|_{Y_1(T)}&\lesssim \varepsilon+(1+T)[a(T)]^p\varepsilon^p+(1+T)^{2p+1}\|u\|_{Y_1(T)}^{pq},\\
	\|v\|_{Y_2(T)}&\lesssim \varepsilon+(1+T)^{2}[a(T)]^{-1}\varepsilon^q+ (1+T)^{q+2}[a(T)]^{pq-1}\|v\|_{Y_2(T)}^{pq},
\end{align*}
where the unexpressed multiplicative constants are independent of $T$ and $\varepsilon$, but may depend on the norms of $u_0,u_1,v_0,v_1$. With the same argument as in Theorem \ref{Thm-Local}, we require that
 	\begin{align*}
 		\begin{cases}
 				(1+T)[a(T)]^p\varepsilon^p\lesssim \varepsilon, \\
 			(1+T)^{2p+1}\varepsilon^{pq}\lesssim\varepsilon, \\ 
 			(1+T)^{2}[a(T)]^{-1}\varepsilon^q\lesssim\varepsilon,\\ 
 			(1+T)^{q+2}[a(T)]^{pq-1}\varepsilon^{pq}\lesssim\varepsilon.
 		\end{cases}
\end{align*}
From these inequalities, we obtain the following upper bound estimates for the existence time $T$:
	\begin{align*}
	T\leqslant 1+T\lesssim\begin{cases}
		\varepsilon^{-\min\left\{\frac{pq-1}{\max\{2p+1,pq+q+1\}},\frac{p-1}{p+1},q-1 \right\}}=\varepsilon^{-\min\left\{\frac{p-1}{p+1},q-1 \right\}}&\mbox{if}\ \ v_1\neq0,\\
		\varepsilon^{-\min\left\{\frac{pq-1}{\max\{2p+1,q+2\}},p-1,\frac{q-1}{2} \right\}}=\varepsilon^{-\min\left\{p-1,\frac{q-1}{2} \right\}}&\mbox{if}\ \ v_1=0,
	\end{cases}
\end{align*}
 	where the unexpressed multiplicative constants are independent of $T$ and $\varepsilon$. Hence, the proof of the lower bound estimates for the lifespan $T(\varepsilon)$ in Theorem \ref{Thm-Local massless case} is complete. 

\begin{remark}
It is possible to prove a local in-time existence result with a higher regularity with respect to time for the $u$-component.

 Indeed, considering
\begin{align*}
u\in Y_1(T):=\ml{C}\big([0,T], H^1_{\ml{L}}(\G)\big)\cap \ml{C}^1\big([0,T], L^2(\G)\big)
\end{align*}
equipped with the norm
\begin{align*}
\|u\|_{Y_1(T)}:=\sup\limits_{t\in[0,T]}\left(\,\sum\limits_{(j,k)\in\Lambda}(1+t)^{\frac{2j+k}{2}}\|\partial_t^j(-\ml{L})^{k/2}u(t,\cdot)\|_{L^2(\G)}\right)
\end{align*}
instead of the norm in \eqref{Norm-Y1-new} and repeating similar computations to the former ones, since $$\int_{t/2}^t(1+t-s)^{-1}\,\mathrm{d}s\lesssim\ln(\mathrm{e}+t),$$ we have to replace \eqref{N1u-v-2} by
\begin{align*}
\|N_1v\|_{Y_1(T)}&\lesssim \varepsilon\|(u_0,u_1)\|_{H^1_{\ml{L}}(\G)\times L^2(\G)} +
(1+T)[a(T)]^p\ln(\mathrm{e}+T)\|v\|_{Y_2(T)}^p.
\end{align*}
From the previous estimate it is clear that we would have a logarithmic loss in lower bound estimates for the lifespan in comparison to those from Theorem \ref{Thm-Local massless case}.
\end{remark}

\section{Concluding remarks}
\hspace{5mm}In the present paper, we proved that local in-time solutions to \eqref{Eq-Main-Systems} blow up in finite time for any $p,q>1$ and under suitable sign assumptions for the Cauchy data, both in the damped Klein-Gordon case and in the massless case. The result in the massless case is expected, due the corresponding blow-up results for the single semilinear damped and undamped wave equations (cf. \cite{Palmieri-JDE=2021, Palmieri-Pseudo=2021}). In the damped Klein-Gordon case, in spite of the fact that for the $L^2(\mathbb{G})$-norm of the solution to the linearized $u$-equation (and of its first order derivatives) we have an exponential decay rate, the coupling with the wave equation destroys the advantageous effect of the mass term.

Furthermore, we also obtained upper and lower estimates for the lifespan of a local in-time solution, depending on the magnitude of the Cauchy data, which is expressed through the positive parameter $\varepsilon$. In the damped Klein-Gordon case (under the further technical assumption $p+q^{-1}<2$ when $v_1=0$), the lifespan estimates that we proved are sharp, in the sense that, up to the multiplicative constant (which is independent of $\varepsilon$), the obtained upper and lower bounds for $T(\varepsilon)$ are the same. In the massless case, since we have less freedom in prescribing a loss of decay for the component $u$ of the solution, our lower bounds for the lifespan do not match the upper bounds, and we cannot claim the sharpness of the lifespan estimates.

Finally, we point out that our blow-up theorems are consistent with those in the Euclidean framework. Indeed, combining the results from  \cite{Chen-Reissig=2021, Wakasugi=2017}, we have that for the Nakao problem with power nonlinearities $(|v|^p,|u|^q)^{\mathrm{T}}$ in $\mathbb{R}^n$, under suitable assumptions for the Cauchy data, the local in-time solutions blow up for any $p,q>1$ satisfying
\begin{align}\label{blow-up range Nakao P in Rn}
\frac{\max\{2+p^{-1},q/2+1\}}{pq-1}>\frac{n-1}{2}.
\end{align} As we mentioned in the introduction, for a compact Lie group $\mathbb{G}$ the global dimension is $0$. Therefore, if we replace $n$ by $0$ in \eqref{blow-up range Nakao P in Rn}, this inequality describing the blow-up range is  satisfied for any $p,q>1$, exactly as in Theorems \ref{Thm-Blow-up} and \ref{Thm-Blow-up massless}.

\section*{Acknowledgments}
W. Chen is supported in part by the National Natural Science Foundation of China (grant No. 12301270), Guangdong Basic and Applied Basic Research Foundation (grant No. 2025A1515010240). A. Palmieri is partially supported by the PRIN 2022 project
``Anomalies in partial differential equations and applications'' CUP H53C24000820006. A. Palmieri
is member of the Gruppo Nazionale per L'Analisi Matematica, la Probabilit\`a e le loro Applicazioni
(GNAMPA) of the {\it Instituto Nazionale di Alta Matematica} (INdAM).

\end{document}